\providecommand{\AddToHook}[2]{#2}
\def\siscbreakurlbefore{\let\siscnativeifpdf\ifpdf\let\ifpdf\iftrue}
\def\siscbreakurlafter{\let\ifpdf\siscnativeifpdf}
\newif\ifsiamclass
\IfFileExists{siamart251216.cls}{\siamclasstrue}{\siamclassfalse}
\ifsiamclass
  \ifdefined\XeTeXversion
    \PassOptionsToClass{nohypdvips}{siamart251216}
    \AddToHook{package/breakurl/before}{\siscbreakurlbefore}
    \AddToHook{package/breakurl/after}{\siscbreakurlafter}
  \fi
  \documentclass[onefignum,onetabnum]{siamart251216}
\else
  \documentclass[11pt]{article}
\fi

\usepackage{amsmath,amssymb,mathtools}
\usepackage{bm}
\usepackage{booktabs}
\usepackage{graphicx}
\usepackage{algorithm}
\usepackage{algpseudocode}

\ifsiamclass
  \newsiamremark{remark}{Remark}
  \newsiamremark{assumption}{Assumption}
\else
  \usepackage{amsthm}
  \usepackage{hyperref}
  \usepackage{cleveref}
  \newtheorem{theorem}{Theorem}
  \newtheorem{lemma}[theorem]{Lemma}
  \newtheorem{proposition}[theorem]{Proposition}
  \newtheorem{corollary}[theorem]{Corollary}
  
  \theoremstyle{remark}
  \newtheorem{remark}[theorem]{Remark}
  \newtheorem{assumption}[theorem]{Assumption}
  \newenvironment{keywords}{\par\noindent\textbf{Keywords. }}{\par}
  \newenvironment{MSCcodes}{\par\noindent\textbf{MSC codes. }}{\par}
  \providecommand{\headers}[2]{}
  
  \providecommand{\email}[1]{\texttt{#1}}

  \let\projectincludegraphics\includegraphics
  \renewcommand{\includegraphics}[2][]{%
    \IfFileExists{#2}{%
      \projectincludegraphics[#1]{#2}%
    }{%
      \fbox{\parbox[c][2in][c]{0.9\linewidth}{\centering
        Figure unavailable in standalone preview:\\[0.5ex]
        \texttt{\detokenize{#2}}}}%
    }%
  }
\fi

\newcommand{\R}{\mathbb R}

\newcommand{\norm}[1]{\left\lVert #1 \right\rVert}

\headers{Stochastic Multilevel Sampling}
        {Zichao Wendy Di}

\title{Stochastic Multilevel Sampling for Inverse Problems via
  Coarse-Space Identifiability}

\author{Zichao Wendy Di\thanks{Mathematics and Computer Science, Argonne National Lab, Lemont, IL
    (\email{wendydi@anl.gov}).}}
\begin{document}
\maketitle
\begin{abstract}
Sampling only part of the data reduces the cost of each iteration in an
inverse problem, but can make reconstruction updates unreliable.
We study stochastic multilevel sampling, which computes updates on a
coarse-to-fine hierarchy of correction spaces.
The guiding idea is that a small sample may support a reliable coarse-grid
update even when it provides insufficient information for a fine-grid update.
For quadratic objectives, we identify conditions under which corrections
computed from sampled data decrease the full-data objective.
In small-batch experiments with computed tomography and periodic deblurring,
the stochastic multilevel method improves reconstruction accuracy over
matched single-level methods.  In an increasingly ill-conditioned CT
refinement with a fixed number of projection angles, the single-level error
increases as the image dimension grows, whereas a deep hierarchy maintains or
improves accuracy and outperforms full-data single-level optimization at the
same number of forward and adjoint applications.  The gains persist across
the tested local solvers and strengthen with deeper hierarchies.
A separate known-probe ptychography study provides complementary evidence for
the benefit of hierarchical corrections for nonlinear inverse problems.
\end{abstract}



\section{Introduction}
\label{sec:introduction}

Large inverse problems are often expressed as sums of many data terms, such as
rays in computed tomography (CT), pixels in an image-deblurring model, or scan
positions in ptychography.  Using only a subset of these terms can substantially
reduce the cost of an iteration~\cite{robbins1951,bottou2018}.  However, the
resulting batch may not contain enough information to determine a stable update
in the full solution space.  This limitation is particularly important for
ill-conditioned inverse problems with weakly observed directions, where a
decrease in the sampled residual need not improve the reconstruction
~\cite{engl1996,hansen2010}.

The reliability of a sampled update depends not only on the number of sampled
data terms, but also on the space in which the update is sought.  A small batch
may be inadequate for a high-dimensional correction while still identifying a
useful correction in a lower-dimensional, structured subspace.  This observation
suggests matching the dimension and geometry of the correction space to the
information available in the sample.  A hierarchy of spaces provides a natural
way to do so, moving from inexpensive, strongly restricted corrections to more
expressive fine-scale corrections as needed.

Classical multigrid provides the best-known example of such a hierarchy: it
combines relaxation with coarse-grid correction to treat error components at
different spatial scales~\cite{brandt1977,xu1992,trottenberg2001}.  Its usual
interpretation, however, concerns the acceleration of a deterministic solve.
Many inverse-problem operators are nonlocal, so the standard connection between
spatial scale and error reduction may not carry over directly
~\cite{natterer2001,kak2001}.  Here we use the hierarchy for a different
purpose.  We ask whether restricting a correction to a coarse space can make a
\emph{sampled} inverse problem statistically more reliable, regardless of
whether the same hierarchy accelerates the corresponding full-data solve.

To formalize this question, let \(H\succeq0\) denote a curvature model at the
current iterate, let \(H_S\) be its sampled counterpart, and let
\(P\in\R^{n\times n_c}\), with \(n_c<n\), define a correction space.  Although
\(H_S\) may approximate \(H\) poorly on \(\R^n\), the restricted matrix
\(P^\top H_S P\) may still approximate \(P^\top H P\) accurately.  We refer to
this property as \emph{coarse-space identifiability}.  The restriction acts on
the unknown variables; it neither coarsens the observations nor changes the
data-acquisition geometry.

In this study, we develop stochastic multilevel sampling around this principle.
The method visits a prescribed hierarchy from coarse to fine, computing each
correction by a gradient update or a short Krylov solve.  A separate data batch
then evaluates the feasible candidates through a projected Armijo line search.
This independent line search limits reuse of the data that constructed the correction,
while the central question remains whether the sampled model identifies a
useful correction on the chosen space.  The sampling distribution, batch size,
and hierarchy remain fixed during a run.

Our analysis follows the correction through three stages.  For quadratic
objectives, an exact identity quantifies the best decrease attainable in a
prescribed correction space and the loss caused by an inexact correction.
Perturbation and concentration bounds then give sufficient conditions under
which the minimizer of the sampled restricted problem decreases the full-data
objective.  Finally, an independent-batch bound gives a sufficient condition
for accepting the projected candidate.

The numerical study first considers an increasingly ill-conditioned CT
refinement in which the image dimension grows while the number of projection
angles remains fixed.  As the measurement-to-unknown ratio decreases, the
single-level sampled method degrades, whereas the multilevel methods maintain
or improve reconstruction accuracy.  Broader sampled-row experiments for CT
and periodic deblurring
then show lower error than matched single-level controls after the same number
of sampled rows or, separately, the same number of forward and adjoint
applications.  The improvement occurs for gradient and short Krylov
corrections and generally increases with hierarchy depth.  A uniform-sampling
deblurring control further shows that nonuniform sampling is not required in
this setting.

We also use known-probe MAGPIE ptychography~\cite{zhang2026magpie} as a
nonlinear scope study, motivated by earlier multigrid optimization for
large-scale ptychographic phase retrieval~\cite{fung2020}.  The hierarchy
produces lower reconstruction error than the single-level method after the
same number of scan passes and, separately, after the same number of updates.
The MAGPIE implementation does not use the sampling strategy or projected
Armijo line search developed here.  These experiments therefore demonstrate
the usefulness of hierarchical correction spaces in a nonlinear problem, not
the performance of the complete stochastic multilevel algorithm.

The remainder of the paper is organized as follows.
\Cref{sec:related} positions the method relative to multilevel optimization and
randomized least-squares methods, and \cref{sec:model-method} defines the
algorithm.  \Cref{sec:theory} develops the coarse-space analysis.
\Cref{sec:experiments} presents the coarse-space diagnostics and
reconstruction experiments, respectively, while \cref{sec:scope} reports the
ptychography study.  \Cref{sec:discussion} summarizes the conclusions and
limitations.

\section{Related methods and positioning}
\label{sec:related}

Multigrid and subspace correction organize updates through hierarchical
spaces~\cite{brandt1977,xu1992}.  MG/OPT and its trust-region and line-search
variants extend this construction to optimization~\cite{nash2000,gratton2008,wen2009}.
In particular, line-search multigrid combines direct and recursive directions
with conditions ensuring descent~\cite{wen2009}.  Our method retains hierarchical
corrections but constructs and evaluates them using separate sampled observations.  Its
additional question is statistical: does a batch represent the restricted
model well enough to produce a useful correction?

Two recent methods combine multilevel structure with stochastic or geometric
optimization.  Marini, Porcelli, and Riccietti~\cite{marini2024multilevel}
use random regularized first-order models and hierarchical descriptions in
variable or function space to obtain an almost-sure global convergence result
for unconstrained optimization.  Vanmaele, Elshiaty, and
Petra~\cite{vanmaele2025multilevel} analyze deterministic geometric coarse
models, including box-constrained tomography, for problems whose objective and
gradient can be evaluated on each spatial grid.  Our setting differs from both:
we retain the same fine-grid inverse-problem objective and acquisition geometry,
sample its measurement rows, and use the spatial hierarchy only to define
correction spaces.  The central issue is therefore whether sampled restricted
gradients and curvature identify a useful correction in a prescribed subspace,
rather than convergence of a stochastic regularized model or deterministic
acceleration through grid-dependent objectives.

Stochastic gradients and variance reduction lower per-step cost and stabilize
updates~\cite{robbins1951,johnson2013,bottou2018}.  Randomized Kaczmarz and
sketch-and-project methods similarly sample rows or
constraints~\cite{strohmer2009,gower2015,needell2016}.  We retain row sampling
but restrict corrections to hierarchical trial spaces.
Stochastic line searches give complementary guarantees under probabilistically
accurate estimates~\cite{mahsereci2017,paquette2020,vaswani2019}; our analysis
connects an independent-batch test of a realized multilevel direction to
fine-objective descent.

Least-squares and Hessian sketching provide the closest theoretical
tools~\cite{woodruff2014,pilanci2016,pilanci2017,drineas2012}.
Concentration on relevant directions is already central to this literature;
for example, iterative Hessian sketch guarantees depend on the geometry of a
tangent cone and retain the full gradient while sketching
curvature~\cite{pilanci2016}.  Here both the restricted gradient and curvature
are sampled, and the prescribed hierarchy changes the trial space across a
sweep.  We apply standard concentration tools~\cite{tropp2015} to these spaces
and study their effect on reconstruction quality at matched sample and operator counts.
The contribution is the multilevel construction and its analysis and numerical
validation, not a new general-purpose embedding bound.  It is also distinct from
multilevel Monte Carlo~\cite{giles2015}, which combines estimators across
discretizations to reduce the variance of a scalar quantity.

The ptychography diagnostic connects to iterative
engines~\cite{maiden2009,maiden2017}, maximum-likelihood
solvers~\cite{thibault2012,odstrcil2018}, and multigrid optimization for
large-scale ptychographic phase retrieval~\cite{fung2020}.  The nonlinear scope
study uses MAGPIE~\cite{zhang2026magpie}; we compare its surrogate hierarchy
with the same surrogate update at single level.

\section{Model and stochastic multilevel method}
\label{sec:model-method}

\subsection{Inverse problems and sampled objectives}
\label{subsec:importance}

We consider a finite-sum inverse problem
\begin{equation}
    \min_{x\in\mathcal C} F(x),
    \qquad F(x)=\sum_{i=1}^m f_i(x),
    \label{eq:finite-sum}
\end{equation}
where \(\mathcal C\subseteq\R^n\) is a nonempty closed convex admissible set,
and each term
corresponds to an observation or observation block.  The analysis and primary
experiments use linear least squares,
\begin{equation}
    F(x)=\frac12\norm{Ax-b}_2^2
        =\frac12\sum_{i=1}^m(a_i^\top x-b_i)^2,
    \label{eq:objective}
\end{equation}
with \(A\in\R^{m\times n}\) and observed data \(b\in\R^m\).  We write
\(\Pi_{\mathcal C}\) for the Euclidean projection onto \(\mathcal C\).
Closed convexity makes this projection single valued and nonexpansive.  The
method does not require a particular constraint: taking \(\mathcal C=\R^n\)
recovers the unconstrained algorithm.  In the CT and deblurring experiments we
use \(\mathcal C=\R_+^n\), since the unknown represents a nonnegative density
or intensity.  This is a modeling choice rather than a property required by the
multilevel construction.  No explicit regularizer is used; inexact corrections,
restricted spaces, and finite iteration counts provide implicit regularization.
For the synthetic linear experiments with ground truth \(x^\star\), we write
\(b=Ax^\star+\xi\), define the relative noise level by
\(\eta=\norm{\xi}_2/\norm{Ax^\star}_2\), and denote the known noise norm by
\(\delta=\norm{\xi}_2\).

Draw a batch \(S=(i_1,\ldots,i_s)\) independently with replacement from a law
\(p_i>0\), where \(s=|S|\), and define
\begin{equation}
    \widetilde F_{S,p}(x)
    =\frac1s\sum_{j=1}^s
       \frac{f_{i_j}(x)}{p_{i_j}},
    \qquad
    \mathbb E[\widetilde F_{S,p}(x)]=F(x).
    \label{eq:importance-objective}
\end{equation}
For least squares, the weighted rows
\begin{equation}
    (\widetilde A_S)_{j,:}
       =\frac{a_{i_j}^\top}{\sqrt{s p_{i_j}}},
    \qquad
    (\widetilde b_S)_j
       =\frac{b_{i_j}}{\sqrt{s p_{i_j}}}
    \label{eq:sample-scaling}
\end{equation}
give \(\widetilde F_{S,p}(x)
=\tfrac12\norm{\widetilde A_Sx-\widetilde b_S}_2^2\), with unbiased gradient
and Hessian.  We use the fixed row-norm law
\(p_i\propto\norm{a_i}_2^2\).  Zero operator rows are omitted from sampling
because they contribute only a constant to the objective.  Uniform
\(p_i=1/m\) gives the usual \(m/s\) rescaling.

\subsection{Hierarchy of correction spaces}

Let \(P_\ell\in\R^{n\times n_\ell}\) prolong level-\(\ell\) corrections to
the finest grid, with \(P_1=I\).  At an iterate \(x\), the restricted full
and sampled objectives are
\begin{equation}
    \phi_\ell(e;x)=F(x+P_\ell e),
    \qquad
    \widetilde\phi_{\ell,S,p}(e;x)=\widetilde F_{S,p}(x+P_\ell e).
    \label{eq:restricted-objectives}
\end{equation}
The hierarchy changes the correction space, not the observation geometry:
ray, pixel, and scan-position indices retain their physical meaning.  In the
CT and deblurring experiments, we transfer from an image grid to the next
coarser grid using the normalized stencil
\begin{equation}
    \frac1{64}
    \begin{bmatrix}
    1 & 3 & 3 & 1\\
    3 & 9 & 9 & 3\\
    3 & 9 & 9 & 3\\
    1 & 3 & 3 & 1
    \end{bmatrix}.
    \label{eq:restriction-stencil}
\end{equation}
In each coordinate, the four weights act at offsets
\(-1,0,1,2\) from the stored coarse-grid anchor, so the stencil is centered at
the intervening half-grid location.  Its interior weights sum to one; entries
outside the image are omitted at the boundary.  Prolongation to the adjacent
finer grid is \(8\) times the transpose of this operation.  The factor \(8\)
is a scaling convention and does not change the range of the prolongation.
We use this choice throughout to define a concrete geometric hierarchy, but
neither the algorithm nor the analysis is restricted to this stencil.  Other
prescribed transfer operators may be used to compose the adjacent
prolongations \(P_\ell\) and construct the correction spaces.

\subsection{Stochastic multilevel cycle}
\label{subsec:kamis-cycle}

A sweep visits levels \(L,L-1,\ldots,1\).  At level \(\ell\), we draw a
\emph{correction batch} \(S\) and use its sampled residual and operator to
compute a local correction using either one gradient step or a short zero-start
LSQR solve~\cite{paige1982}.  Prolongation maps this correction to a fine-space
direction \(d\).  The correction batch also supplies the scale
\begin{equation}
  \bar\alpha_S
  =\max\!\left\{0,
    \frac{r_S^\top \widetilde A_Sd}{\norm{\widetilde A_Sd}_2^2}
  \right\},
  \qquad r_S=\widetilde b_S-\widetilde A_Sx,
  \label{eq:correction-batch-scale}
\end{equation}
with \(\bar\alpha_S=0\) if the denominator vanishes.  This value is used only
to initialize a finite backtracking grid,
\begin{equation}
  \mathcal B_S=\{0,\bar\alpha_S,\bar\alpha_S/2,\ldots,
    \bar\alpha_S/2^J\}.
  \label{eq:projected-candidate-set}
\end{equation}
We use \(J=5\).  The construction is invariant to positive rescaling of \(d\)
and fixes \(\mathcal B_S\) before validation data are drawn.

We then draw an independent \emph{line-search batch} \(C\) and evaluate the
projected path
\begin{equation}
  z(\beta)=\Pi_{\mathcal C}(x+\beta d),
  \qquad \beta\in\mathcal B_S.
  \label{eq:projected-search-path}
\end{equation}
Starting from \(\bar\alpha_S\), we accept the first candidate that satisfies a
projected Armijo condition on the line-search batch.  If no candidate does, the
iterate is unchanged.  Thus projection is part of the search itself, rather
than a modification applied after an unconstrained step has been selected.
Degenerate directions leave the iterate unchanged.
\Cref{subsec:stochastic-line-search} analyzes the projected-candidate test.

\subsection{Complete procedure}

\Cref{alg:amls} collects the complete procedure for a prescribed number of
coarse-to-fine sweeps.  The numerical experiments state their endpoint
conventions separately.

\begin{algorithm}[htbp]
\caption{Stochastic multilevel sampling with a projected Armijo line search.}
\label{alg:amls}
\begin{algorithmic}[1]
\Require \(A,b,\mathcal C,\{P_\ell\}_{\ell=1}^L,p,s_S,s_C,x^{(0)}\);
  local solver (gradient or LSQR(\(t\))), Armijo constant \(c\in(0,1)\),
  backtracking limit \(J\), and number of sweeps \(K\)
\State \(x\gets\Pi_{\mathcal C}(x^{(0)})\)
\For{\(k=1,2,\ldots,K\)}
  \For{\(\ell=L,L-1,\ldots,1\)}
    \State draw \(s_S\) correction rows \(S\sim p\); form
      \(\widetilde A_S,\widetilde b_S\) using \cref{eq:sample-scaling}
    \State \(G\gets\widetilde A_SP_\ell\),
      \(r\gets\widetilde b_S-\widetilde A_Sx\)
    \State compute \(e=G^\top r\), or \(e\) by \(t\) zero-start LSQR steps
      for \(\min_z\norm{Gz-r}_2\); set \(d\gets P_\ell e\)
    \State \(v_S\gets\widetilde A_Sd\)
    \If{\(d=0\) or \(v_S=0\)}
      \State \(\bar\alpha_S\gets0\)
    \Else
      \State \(\bar\alpha_S\gets
        \max\{0,r^\top v_S/\norm{v_S}_2^2\}\)
    \EndIf
    \State \(\mathcal B_S\gets
      \{0,\bar\alpha_S,\bar\alpha_S/2,\ldots,\bar\alpha_S/2^J\}\)
    \State draw \(s_C\) independent line-search rows \(C\sim p\); form
      \(\widetilde A_C,\widetilde b_C\) using \cref{eq:sample-scaling}
    \State initialize \(\beta_\star\gets0\), \(x_{\rm trial}\gets x\)
    \For{\(\beta=\bar\alpha_S,\bar\alpha_S/2,\ldots,
      \bar\alpha_S/2^J\)}
      \State \(z_\beta\gets\Pi_{\mathcal C}(x+\beta d)\)
      \If{\(\nabla\widetilde F_{C,p}(x)^\top(z_\beta-x)<0\) and
        \(\widetilde F_{C,p}(z_\beta)-\widetilde F_{C,p}(x)
        \leq c\nabla\widetilde F_{C,p}(x)^\top(z_\beta-x)\)}
        \State \(\beta_\star\gets\beta\), \(x_{\rm trial}\gets z_\beta\);
          \textbf{break}
      \EndIf
    \EndFor
    \State \(x\gets x_{\rm trial}\)
  \EndFor
\EndFor
\State \Return \(x\)
\end{algorithmic}
\end{algorithm}

\section{Coarse-space mechanism}
\label{sec:theory}

We analyze when a sampled coarse correction decreases a quadratic objective,
without tying the mechanism to a particular local solver.
Objective decrease is separate from reconstruction accuracy in noisy inverse
problems; the experiments assess the latter.  The results also apply to local
quadratic models, but do not establish convergence for a nonlinear problem.

\subsection{Exact and sampled coarse-space decrease}

Let
\begin{equation}
    F(x)=\frac12 x^\top Hx-h^\top x,
    \qquad
    H\succeq0,
    \label{eq:quadratic-model}
\end{equation}
and let \(P\in\R^{n\times n_c}\) be a coarse-space basis.  We define
\begin{equation}
    B=P^\top H P,
    \qquad
    B_S=P^\top H_S P,
    \qquad
    c=P^\top g(x),
    \qquad
    c_S=P^\top g_S(x),
    \label{eq:projected-quantities}
\end{equation}
where \(H_S\) and \(g_S\) are the sampled Hessian and sampled gradient.

\begin{lemma}[Exact coarse-space decrease]
\label{lem:exact-coarse-decrease}
For the quadratic model \cref{eq:quadratic-model}, fix \(x\) and a coarse basis
\(P\) such that \(B=P^\top H P\) is positive definite.  Let
\(c=P^\top\nabla F(x)\) and \(e_\star=-B^{-1}c\).  Then every coarse correction
\(e\in\R^{n_c}\) satisfies
\begin{equation}
    F(x)-F(x+Pe)
    =
    -c^\top e-\frac12 e^\top B e
    =
    \frac12\norm{c}_{B^{-1}}^2
    -
    \frac12\norm{e-e_\star}_B^2 .
    \label{eq:exact-coarse-decrease}
\end{equation}
Consequently, when \(c\neq0\), define the normalized coarse-decrease capture by
\begin{equation}
    \kappa(e;P,x)
    :=
    \frac{F(x)-F(x+Pe)}{\tfrac12\norm{c}_{B^{-1}}^2}
    =
    1-\frac{\norm{e-e_\star}_B^2}{\norm{e_\star}_B^2}
    \leq 1,
    \label{eq:exact-capture-ratio}
\end{equation}
which may be negative.  The correction is a descent step if and only if
\(\kappa(e;P,x)>0\), equivalently
\(\norm{e-e_\star}_B<\norm{e_\star}_B\).
\end{lemma}

\begin{proof}
Expanding \(F(x+Pe)-F(x)=c^\top e+\tfrac12e^\top Be\) and completing the
square with \(e_\star=-B^{-1}c\) gives \cref{eq:exact-coarse-decrease}.
Division by \(\tfrac12\norm{c}_{B^{-1}}^2
=\tfrac12\norm{e_\star}_B^2\) gives \cref{eq:exact-capture-ratio}.
\end{proof}

The sampled Hessian and gradient are formed from the same rows, so their
errors should not generally be treated as independent perturbations.  The
following result retains their coupling.

\begin{theorem}[Coupled sampled-to-full transfer]
\label{thm:coupled-sampled-to-full-transfer}
Suppose that \(B\) and \(B_S\) are positive definite, and let
\(e_\star=-B^{-1}c\) and \(e_S=-B_S^{-1}c_S\).  Define the coupled
normal-equation defect and sampled-Hessian stability factor by
\begin{equation}
\begin{aligned}
    r_S &:= (c_S-c)+(B_S-B)e_\star,\\
    \zeta_S &:= \norm{r_S}_{B^{-1}},
    &\omega_S&:=\norm{B^{1/2}B_S^{-1}B^{1/2}}_2.
\end{aligned}
\label{eq:coupled-defect}
\end{equation}
Then
\begin{equation}
    F(x)-F(x+Pe_S)
    \geq
    \frac12\norm{c}_{B^{-1}}^2
    -\frac12\omega_S^2\zeta_S^2.
    \label{eq:coupled-transfer-bound}
\end{equation}
Consequently, \(\omega_S\zeta_S<\norm{c}_{B^{-1}}\) is sufficient for
the sampled correction to decrease the full objective.
\end{theorem}

\begin{proof}
The two normal equations give the exact relation
\[
    e_S-e_\star=-B_S^{-1}r_S.
\]
Therefore,
\[
    \norm{e_S-e_\star}_B
    \leq
    \norm{B^{1/2}B_S^{-1}B^{1/2}}_2
    \norm{r_S}_{B^{-1}}
    =\omega_S\zeta_S.
\]
Substitution into \cref{eq:exact-coarse-decrease} proves
\cref{eq:coupled-transfer-bound}.
\end{proof}

For \(c\ne0\), the corresponding normalized certificate is
\begin{equation}
    \kappa(e_S;P,x)
    \geq
    1-\left(
      \frac{\omega_S\zeta_S}{\norm{c}_{B^{-1}}}
    \right)^2.
    \label{eq:coupled-capture-bound}
\end{equation}
Here \(\zeta_S\) measures how poorly the ideal full-data correction satisfies
the sampled normal equations, while \(\omega_S\) measures how strongly the
sampled solve can amplify that defect.  Their product bounds the distance
between the sampled and ideal corrections.

\begin{assumption}[Coarse-space concentration]
\label{as:concentration}
The projected Hessian \(B\) is positive definite.  There is a constant
\(0\leq\delta<1\) such that the sampled projected Hessian obeys the relative
spectral bound
\begin{equation}
    (1-\delta)B\preceq B_S\preceq(1+\delta)B,
    \label{eq:hessian-concentration}
\end{equation}
and there is a constant \(\varepsilon\geq0\) such that the sampled projected
gradient obeys
\begin{equation}
    \norm{c_S-c}_{B^{-1}}\leq\varepsilon .
    \label{eq:gradient-concentration}
\end{equation}
\end{assumption}

This assumption does not require \(H\) to be positive definite on the full
space; it excludes null directions only from the chosen correction space.
For \(P=P_\ell\), write
\(B_\ell=P_\ell^\top H P_\ell\) and
\(B_{\ell,S}=P_\ell^\top H_S P_\ell\).  The corresponding relative
sampled-Hessian error is
\begin{equation}
    \delta_\ell(S)
    :=
    \norm{
      B_\ell^{-1/2}(B_{\ell,S}-B_\ell)B_\ell^{-1/2}
    }_2.
    \label{eq:relative-hessian-error}
\end{equation}
It is the smallest value of \(\delta\) for which
\cref{eq:hessian-concentration} holds on level \(\ell\).

\begin{corollary}[Separate-error transfer bound]
\label{thm:sampled-to-full-transfer}
Under \cref{as:concentration}, the full objective decrease obtained by the
sampled coarse correction \(e_S=-B_S^{-1}c_S\) satisfies
\begin{equation}
    F(x)-F(x+Pe_S)
    \geq
    \frac12\norm{c}_{B^{-1}}^2
    -
    \frac12
    \left(
    \frac{\delta}{1-\delta}\norm{c}_{B^{-1}}
    +
    \frac{\varepsilon}{1-\delta}
    \right)^2 .
    \label{eq:sampled-to-full-transfer-bound}
\end{equation}
\end{corollary}

\begin{proof}
The concentration condition \cref{eq:hessian-concentration} implies
\[
    \omega_S\leq\frac{1}{1-\delta}.
\]
Moreover,
\[
    \zeta_S
    \leq
    \norm{c_S-c}_{B^{-1}}
    +\norm{(B_S-B)e_\star}_{B^{-1}}
    \leq
    \varepsilon+\delta\norm{e_\star}_B
    =\varepsilon+\delta\norm{c}_{B^{-1}}.
\]
Substitution into the coupled bound \cref{eq:coupled-transfer-bound} gives
\cref{eq:sampled-to-full-transfer-bound}.
\end{proof}

For \(c\ne0\), \cref{thm:sampled-to-full-transfer} gives
\[
    \kappa(e_S;P,x)
    \geq
    1-\left(
      \frac{\delta}{1-\delta}
      +
      \frac{\varepsilon}{(1-\delta)\norm{c}_{B^{-1}}}
    \right)^2.
\]
This lower bound is positive precisely when
\(2\delta+\varepsilon/\norm{c}_{B^{-1}}<1\).
Thus \(\delta<1\) ensures that the sampled coarse Hessian is invertible, but
does not by itself ensure descent.  A nonpositive bound also does not imply
that a realized correction fails.
The separate-error bound follows from a triangle inequality applied to the
coupled defect.  It can therefore be much more conservative than
\cref{eq:coupled-capture-bound}, particularly when the Hessian and gradient
errors formed from the same sampled rows cancel in the normal equations.
For \(P=P_\ell\), we write the realized level-\(\ell\) capture as
\(\kappa_\ell(S):=\kappa(e_S;P_\ell,x)\).

\subsection{Sample concentration on the hierarchy}

For row-sampled least squares, let \(v_i=P_\ell^\top a_i\) and
\(s_i=v_i^\top B_\ell^{-1}v_i\) be the coarse row and its leverage score,
respectively.  Since \(\sum_i s_i=n_\ell\), the coarse-space
coherence is~\cite{drineas2012}
\begin{equation}
    \mu_\ell = \frac{m}{n_\ell}\max_{1\le i\le m} s_i \ge 1 .
    \label{eq:coarse-coherence}
\end{equation}
For a sampling law \(p_i>0\), define the coarse importance mismatch
\begin{equation}
    \Lambda_\ell(p)=\max_{1\leq i\leq m}\frac{s_i}{p_i}.
    \label{eq:importance-mismatch}
\end{equation}
Uniform sampling gives \(\Lambda_\ell(p)=\mu_\ell n_\ell\), while the ideal
coarse-leverage law \(p_i=s_i/n_\ell\) on nonzero-leverage rows gives
\(\Lambda_\ell(p)=n_\ell\).

\begin{proposition}[Coarse-space concentration]
\label{prop:concentration}
Fix a level \(\ell\), and draw \(s\) rows independently with replacement from
a law \(p_i>0\).  Form the unbiased reweighted sampled Hessian
\[
    B_{\ell,S,p}
    = \frac1s\sum_{j=1}^s
      \frac{1}{p_{i_j}}v_{i_j}v_{i_j}^\top .
\]
For any \(\delta\in(0,1)\) and \(\rho\in(0,1)\), if
\begin{equation}
    s \;\ge\;
        \frac{3\,\Lambda_\ell(p)}{\delta^2}\,
        \log\frac{2 n_\ell}{\rho},
    \label{eq:sample-window}
\end{equation}
then the Hessian bound \cref{eq:hessian-concentration} of
\cref{as:concentration}, with \(B_S=B_{\ell,S,p}\), holds with probability at
least \(1-\rho\).  For uniform iid sampling,
\(\Lambda_\ell=\mu_\ell n_\ell\).  The diagnostics below use this
with-replacement importance sampler.
\end{proposition}

\begin{proof}[Proof sketch]
Write \(u_i = B_\ell^{-1/2} v_i\), so that \(\sum_i u_i u_i^\top = I_{n_\ell}\)
and
\[
    B_\ell^{-1/2}B_{\ell,S,p}B_\ell^{-1/2}
    = \frac1s\sum_{j=1}^s
      \frac{1}{p_{i_j}}u_{i_j}u_{i_j}^\top .
\]
For one random draw, the matrix \(p_{i}^{-1}u_i u_i^\top\) has expectation
\(I_{n_\ell}\) and spectral norm at most \(\Lambda_\ell(p)\).  Matrix
Chernoff~\cite{tropp2015} therefore gives
\[
    \Pr\!\left(
      \left\|B_\ell^{-1/2}B_{\ell,S,p}B_\ell^{-1/2}-I_{n_\ell}\right\|_2
      >\delta
    \right)
    \le 2n_\ell\exp\!\left(
      -\frac{s\delta^2}{3\Lambda_\ell(p)}
    \right).
\]
The right-hand side is at most \(\rho\) under
\cref{eq:sample-window}.
\end{proof}

\begin{proposition}[Projected-gradient concentration]
\label{prop:projected-gradient-concentration}
Fix a level \(\ell\) and an iterate \(x\).  Let
\(\xi_i(x)=a_i^\top x-b_i\) be the signed residual in the gradient convention,
so that
\[
    c_\ell(x)=P_\ell^\top\nabla F(x)=\sum_{i=1}^m v_i\xi_i(x).
\]
For the same with-replacement importance sampler as in
\cref{prop:concentration}, define the unbiased projected-gradient estimator
\[
    c_{\ell,S,p}(x)
    =
    \frac1s\sum_{j=1}^s
    \frac{1}{p_{i_j}}v_{i_j}\xi_{i_j}(x).
\]
Then
\begin{equation}
    \mathbb E\!\left[
      \norm{c_{\ell,S,p}(x)-c_\ell(x)}_{B_\ell^{-1}}^2
    \right]
    =
    \frac{1}{s}
    \left(
      \sum_{i=1}^m \frac{\xi_i(x)^2s_i}{p_i}
      -
      \norm{c_\ell(x)}_{B_\ell^{-1}}^2
    \right).
    \label{eq:projected-gradient-second-moment}
\end{equation}
Consequently, with probability at least \(1-\rho\),
\begin{equation}
    \norm{c_{\ell,S,p}(x)-c_\ell(x)}_{B_\ell^{-1}}
    \leq
    \left[
      \frac{1}{s\rho}
      \left(
        \sum_{i=1}^m \frac{\xi_i(x)^2s_i}{p_i}
        -
        \norm{c_\ell(x)}_{B_\ell^{-1}}^2
      \right)
    \right]^{1/2}.
    \label{eq:projected-gradient-concentration}
\end{equation}
\end{proposition}

\begin{proof}
Let \(Y=B_\ell^{-1/2}v_I\xi_I(x)/p_I\), where
\(\mathbb P(I=i)=p_i\).  Then
\(\mathbb E Y=B_\ell^{-1/2}c_\ell(x)\), and
\[
    \mathbb E\norm{Y}_2^2
    =
    \sum_{i=1}^m p_i
    \frac{\xi_i(x)^2}{p_i^2}
    v_i^\top B_\ell^{-1}v_i
    =
    \sum_{i=1}^m \frac{\xi_i(x)^2s_i}{p_i}.
\]
For an average of \(s\) independent copies of \(Y\), the mean-square error is
the single-sample variance divided by \(s\), which gives
\cref{eq:projected-gradient-second-moment}.  Markov's inequality gives
\cref{eq:projected-gradient-concentration}.
\end{proof}

\begin{remark}[Inexact local solves]
\label{rmk:krylov-identifiability}
\Cref{thm:coupled-sampled-to-full-transfer} concerns the exact minimizer of the sampled
coarse problem.  The gradient and short-LSQR corrections in
\cref{alg:amls} are inexact, but their realized quality is still measured by
\cref{eq:exact-capture-ratio}.  A short LSQR solve uses a Krylov space built
from the same correction batch, so a concentration result for a fixed or
independently chosen subspace does not apply directly.  We therefore assess
these corrections empirically.
\end{remark}

\subsection{Projected Armijo line search on an independent batch}
\label{subsec:stochastic-line-search}

The correction batch determines a direction, whereas the line-search batch
tests the resulting projected candidates against independent data.  This
separation is needed because correction-batch decrease does not imply
full-objective decrease.  Condition on the iterate \(x\), correction batch
\(S\), direction \(d\), and hence the finite set \(\mathcal B_S\) in
\cref{eq:projected-candidate-set}.  For each \(\beta\in\mathcal B_S\), define
\begin{equation}
  z_\beta=\Pi_{\mathcal C}(x+\beta d),
  \qquad s_\beta=z_\beta-x.
  \label{eq:projected-displacement}
\end{equation}
Although \(s_\beta\) need not be collinear with \(d\), it is fixed before the
line-search batch is drawn.  In particular, projection can move the realized
displacement outside \(\operatorname{range}(P_\ell)\).  The coarse-space theory
therefore explains how the proposal is constructed, while the line-search test
below controls the actual feasible displacement.  Draw
\(C=(i_1,\ldots,i_{s_C})\) with replacement from \(p_i>0\), as in stochastic
line-search methods~\cite{paquette2020,vaswani2019}.  Writing
\(r_i(x)=b_i-a_i^\top x\), let
\begin{equation}
\begin{aligned}
  q_\beta &=r(x)^\top As_\beta,
  &\widehat q_{\beta,C}
    &=\frac1{s_C}\sum_{j=1}^{s_C}
      \frac{r_{i_j}(x)a_{i_j}^\top s_\beta}{p_{i_j}},\\
  h_\beta &=\norm{As_\beta}_2^2,
  &\widehat h_{\beta,C}
    &=\frac1{s_C}\sum_{j=1}^{s_C}
      \frac{(a_{i_j}^\top s_\beta)^2}{p_{i_j}}.
\end{aligned}
\label{eq:projected-linesearch-qh}
\end{equation}
These estimators are conditionally unbiased, and the exact objective change is
\begin{equation}
  F(z_\beta)-F(x)=-q_\beta+\frac12h_\beta.
  \label{eq:projected-objective-change}
\end{equation}

\begin{proposition}[Sufficient decrease for a projected candidate]
\label{prop:projected-candidate-decrease}
Suppose that, simultaneously for every \(\beta\in\mathcal B_S\),
\begin{equation}
  q_\beta\geq\widehat q_{\beta,C}-\varepsilon_{q,\beta},
  \qquad
  h_\beta\leq\widehat h_{\beta,C}+\varepsilon_{h,\beta}
  \label{eq:projected-one-sided-bounds}
\end{equation}
with conditional probability at least \(1-\rho\).  Then any candidate satisfying
\begin{equation}
  -\widehat q_{\beta,C}+\varepsilon_{q,\beta}
  +\frac12\bigl(\widehat h_{\beta,C}+\varepsilon_{h,\beta}\bigr)<0
  \label{eq:projected-confidence-test}
\end{equation}
strictly decreases the full objective with the same conditional probability.
\end{proposition}

\begin{proof}
Substitute \cref{eq:projected-one-sided-bounds} into
\cref{eq:projected-objective-change}.  The bounds are simultaneous, so the
conclusion remains valid when a candidate is selected from the finite set.
\end{proof}

Scalar concentration and a union bound over \(\mathcal B_S\) can provide the
radii in \cref{eq:projected-one-sided-bounds}.  In the experiments, starting from
\(\bar\alpha_S\), they accept the first projected candidate for which
\begin{equation}
  \widehat q_{\beta,C}>0,
  \qquad
  -\widehat q_{\beta,C}+\frac12\widehat h_{\beta,C}
  \leq-c\widehat q_{\beta,C},
  \label{eq:implemented-line-search}
\end{equation}
with \(c=10^{-4}\), and halve the step otherwise.  This is precisely
\(\widetilde F_{C,p}(z_\beta)-\widetilde F_{C,p}(x)
\leq c\nabla\widetilde F_{C,p}(x)^\top(z_\beta-x)<0\).
If no candidate satisfies the condition, rejection leaves the iterate
unchanged.  The rule evaluates the actual feasible update without reusing the
batch that constructed the direction.  It guarantees decrease of the sampled
line-search objective; \cref{prop:projected-candidate-decrease} gives the
additional margin sufficient for decrease of the full objective.  When projection
is inactive, the rule reduces to an ordinary Armijo line search along
\(d\).

\section{Numerical results}
\label{sec:experiments}

The experiments address two questions.  First, can a low-dimensional
correction space be identified more reliably from a small sample than the fine
space?  Second, does that advantage translate into lower reconstruction error
under the current projected Armijo line search?  We study the first question
through correction diagnostics and sufficient-decrease bounds, and the second
through matched large-grid reconstruction experiments.

\subsection{Sample identifiability and descent bounds}
\label{subsec:mechanism}

At \(x=0\), we measure curvature error, direction alignment, and realized
decrease for sampled-subproblem solves.  For each CT grid, we evaluate all
correction spaces using the same iid row-norm batch and the
inverse-probability scaling of \cref{subsec:importance}.  Because
small batches produce rank-deficient sampled Hessians, we use the
ridge-stabilized diagnostic objective
\begin{equation}
 f_\ell^\gamma(e)
 =\tfrac12\|AP_\ell e-b\|^2+\tfrac{\gamma_\ell}{2}\|e\|^2,
 \qquad
 \gamma_\ell=10^{-8}\max\{1,\|P_\ell^\top A^\top AP_\ell\|_2\}.
 \label{eq:diagnostic-ridge}
\end{equation}
We add \(\gamma_\ell I\) to both the full and sampled Hessians and report the
resulting curvature error \(\delta_\ell^\gamma(S)\) and capture ratio
\(\kappa_\ell^\gamma(S)\).  This ridge is used only for the diagnostics, not in
the reconstruction algorithm.

\begin{figure}[htbp]
\centering
\includegraphics[width=\textwidth]{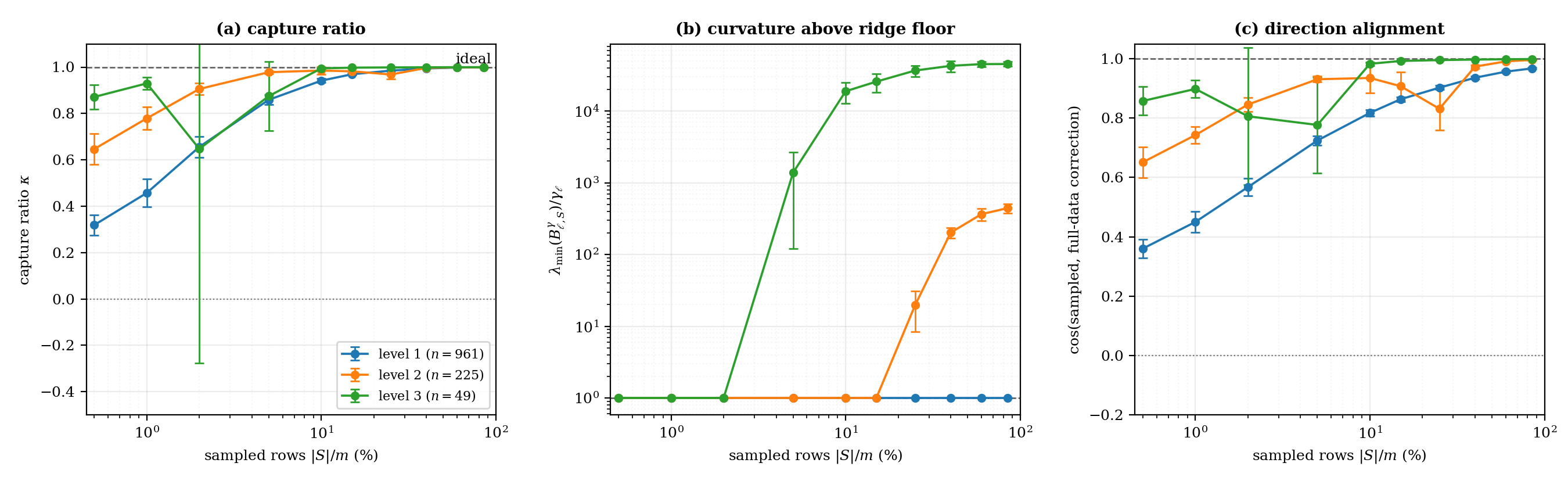}
\caption{Fine-versus-coarse correction diagnostics for \(N=31\) and \(30\)
projection angles, using identical iid importance-weighted batches at every
level.  Markers are means over twelve seeds, and bars show one standard
deviation.  (a) Ridge-model capture ratio; (b) minimum sampled curvature
relative to the ridge floor; (c) cosine between sampled and full-data ridge
corrections.}
\label{fig:fine-vs-coarse}
\end{figure}

\Cref{fig:fine-vs-coarse} compares corrections under the same batch on an
\(N=31\), \(30\)-angle problem with \(m=1320\).  At a \(0.5\%\) batch
\(\lvert S\rvert=7\), the
mean capture ratios on levels 1--3 are \(0.319\), \(0.647\), and \(0.871\);
the fine and coarsest direction cosines are \(0.360\) and \(0.858\).
Thus the coarsest space gives a better aligned correction in this regime.
At larger batches, the coarse sampled curvature rises above the ridge floor and their
corrections approach the corresponding full-data corrections.  Capture and
direction alignment measure the quality of a particular correction; establishing
the entire coarse model additionally requires the concentration conditions.

We next test whether the descent certificates remain informative on larger
problems.  We evaluate full-angle CT grids with \(N\in\{255,511\}\), \(N\)
angles, and fixed \(1\%\) noise.  The two deepest geometric correction spaces
have dimensions \(961\) and \(225\), and we test batch fractions \(0.5\%\),
\(5\%\), and \(50\%\) using six shared sampling seeds.

\Cref{tab:large-concentration} reports how often the coupled lower bound is
positive.  A positive bound certifies that the sampled correction decreases
the full-data ridge objective.  At \(N=255\), the \(0.5\%\) batch contains
\(460\) row draws.  It certifies the \(225\)-dimensional space in five of six
trials but never certifies the \(961\)-dimensional space, which has more
unknown correction coefficients than sampled rows.  Batches of \(5\%\) or
larger certify both spaces in every trial.  At \(N=511\), even the \(0.5\%\)
batch contains \(1847\) draws and certifies both spaces in all six trials.
Thus the experiment directly shows that a smaller correction space can be
identified from fewer sampled rows.  These ridge-quadratic diagnostics are
evaluated only at the initial iterate and are not reconstruction results.

\begin{table}[htbp]
{\footnotesize
\caption{Empirical verification of the coupled descent bound
\cref{eq:coupled-capture-bound}.  Entries give the number of positive bounds
among six trials, evaluated at \(x=0\) and \(\eta=0.01\) using
\eqref{eq:diagnostic-ridge}.  Batch fractions count iid row draws, including
repeats.}
\label{tab:large-concentration}
\begin{center}
\begin{tabular}{rrrr}
\toprule
\(N\) & batch & \(n_\ell=225\) & \(n_\ell=961\) \\
\midrule
255 & \(0.5\%\) & \(5/6\) & \(0/6\) \\
255 & \(5\%\)   & \(6/6\) & \(6/6\) \\
255 & \(50\%\)  & \(6/6\) & \(6/6\) \\
511 & \(0.5\%\) & \(6/6\) & \(6/6\) \\
511 & \(5\%\)   & \(6/6\) & \(6/6\) \\
511 & \(50\%\)  & \(6/6\) & \(6/6\) \\
\bottomrule
\end{tabular}
\end{center}
}
\end{table}

\subsection{Dependence on the correction space}
\label{subsec:subspace}

We also compare the geometric space with same-dimensional DCT low-pass,
Haar-random, leading right-singular-vector, and full-data Krylov spaces; the
complete definitions and results are given in
\cref{app:subspace-comparison}.  For CT, the geometric and DCT spaces have
approximation defects \(0.202\) and \(0.206\) and capture ratios \(0.975\) and
\(0.976\), respectively.  The random space has approximation defect \(0.970\)
and capture \(0.646\); the deblurring comparison shows the same qualitative
separation.

Thus reliable sampling does not compensate for a correction space that poorly
represents the solution.  The comparison supports balancing approximation and
sampling reliability, rather than assigning a universal advantage to
geometric coarsening.  These diagnostics use particular test images at
\(x=0\), and the right-singular-vector and Krylov spaces are references rather
than cost-matched algorithms.

\subsection{Reconstruction with a projected Armijo line search}
\label{sec:kamis}
\label{subsec:kernel-comparison}

We first test whether the hierarchy remains effective as fixed-angle CT
becomes increasingly ill-conditioned under grid refinement, and then
test whether the same behavior persists across CT and deblurring, two grid
sizes, two noise levels, three local solvers, and four batch fractions.  The matched
single-level control
sets \(L=1\) while retaining the same local solver, row-norm sampler,
projected Armijo line search, and nonnegativity constraint.  It is not a tuned
full-data or hybrid-regularized Krylov solver, and no explicit regularizer is
used.

We test three local solvers: gradient, LSQR(2), and LSQR(5).
Each solver is run at single level, with the three-level geometric hierarchy,
and with a deeper hierarchy ending on a \(15\times15\) grid.
The latter has five levels at \(N=255\) and six at \(N=511\).
All variants retain the same row-norm sampling law and projected
Armijo line search on an independent batch.  The correction batch sets the step scale in
\cref{eq:correction-batch-scale}; the line-search batch compares the resulting
projected nonnegative candidates.  A nonzero one-step LSQR correction is collinear with the
sampled-gradient direction, and scaling the candidate grid by
\(\bar\alpha_S\) removes this scalar difference in exact arithmetic.  The
gradient control therefore tests whether hierarchy gains require the
additional local Krylov steps.

To test the method beyond the fixed-angle CT sequence, we perform a matched
sweep at \(N\in\{255,511\}\) for both CT and periodic deblurring.  The CT
problems use \(N\) projection angles.  The deblurring problems use periodic
convolution with an isotropic Gaussian point-spread function whose standard
deviation is \(\sigma=2\) pixels.  Thus the blur width is fixed in pixel units
as the grid grows; this is not a fixed-physical-width refinement.  For both
problems, we use relative noise levels \(\eta\in\{0.01,0.05\}\).

We test iid importance-weighted batch fractions of \(0.5\%\), \(2\%\),
\(5\%\), and \(10\%\).  The \(0.5\%\) fraction is a small-batch stress test,
not an optimized reconstruction parameter.  For each configuration, six
paired replicates redraw the measurement noise and sampled rows.  Within a
replicate, the noisy data are shared across local solvers and hierarchy
choices, and the correction and line-search draws are paired by update index.
The sweep contains \(432\) trajectories per batch fraction and \(1728\)
trajectories in total.

To compare methods at matched computational effort, we use two endpoints.  A
row epoch is \(m\) sampled rows, including the correction and line-search
batches, and an operator epoch is \(m\) scalar-row applications of the forward
operator or its adjoint.  We report the last completed update within twelve
row epochs and, separately, within sixty operator epochs.  The operator count
includes inner LSQR applications, correction-batch residual and scale
evaluations, and projected line-search evaluations; coarse-level applications
receive no discount.  We exclude setup, transfers, and diagnostic evaluations,
so these endpoints normalize the principal operator applications rather than
wall-clock time.

\begin{figure}[htbp]
\centering
\includegraphics[width=\textwidth]{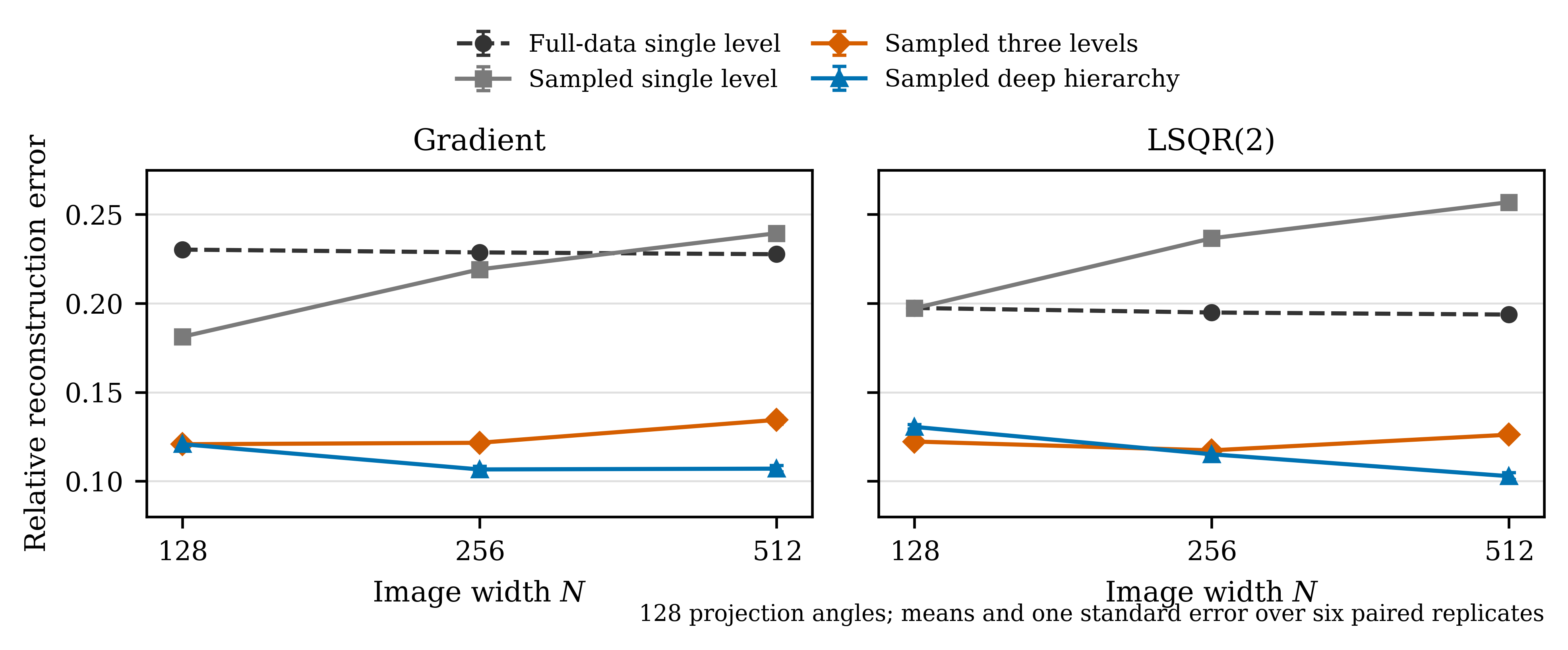}
\caption{Relative reconstruction error under CT grid refinement with \(128\)
fixed projection angles.  All cases use \(5\%\) relative noise, and the
sampled methods use \(5\%\) iid row batches.  Points are means over six paired
replicates at the last complete update within twenty-four normalized forward
and adjoint applications; bars show one standard error.}
\label{fig:fixed-angle-ct}
\end{figure}

\paragraph{Increasingly ill-conditioned CT refinement.}
We construct an increasingly ill-conditioned sequence by holding the number
of CT projection angles at \(128\) as the image width increases through
\(N\in\{128,256,512\}\).  The detector grid contains \(182\), \(363\), and
\(725\) bins per angle, respectively, so the measurement-to-unknown ratios are
\(m/n=1.42\), \(0.71\), and \(0.35\).  Thus the number of measurements grows
fourfold while the number of unknowns grows sixteenfold, and the two finer
systems are underdetermined.  We use \(5\%\) row batches, \(5\%\) relative
noise, six paired replicates, and either gradient or LSQR(2) corrections.
\Cref{fig:fixed-angle-ct} reports the last complete update within twenty-four
normalized forward and adjoint applications.

From \(N=128\) to \(N=512\), the sampled single-level error increases by
\(32\%\) for gradient corrections and \(30\%\) for LSQR(2).  Over the same
refinement, the deep-hierarchy error decreases by \(12\%\) and \(21\%\),
respectively, while the three-level error remains between \(0.117\) and
\(0.134\).  At \(N=512\), the deep hierarchy gives mean errors \(0.1070\) and
\(0.1028\) for gradient and LSQR(2), compared with \(0.2276\) and \(0.1937\)
for the corresponding full-data single-level controls.  These reductions
occur in all six paired replicates.  The three grid sizes do not establish
grid-independent convergence, but they show that the hierarchy prevents the
loss of accuracy observed for the sampled single-level method as angular
information becomes increasingly sparse relative to the image dimension.  The
comparison counts operator applications rather than wall-clock time.

\begin{figure}[htbp]
\centering
\includegraphics[width=\textwidth]{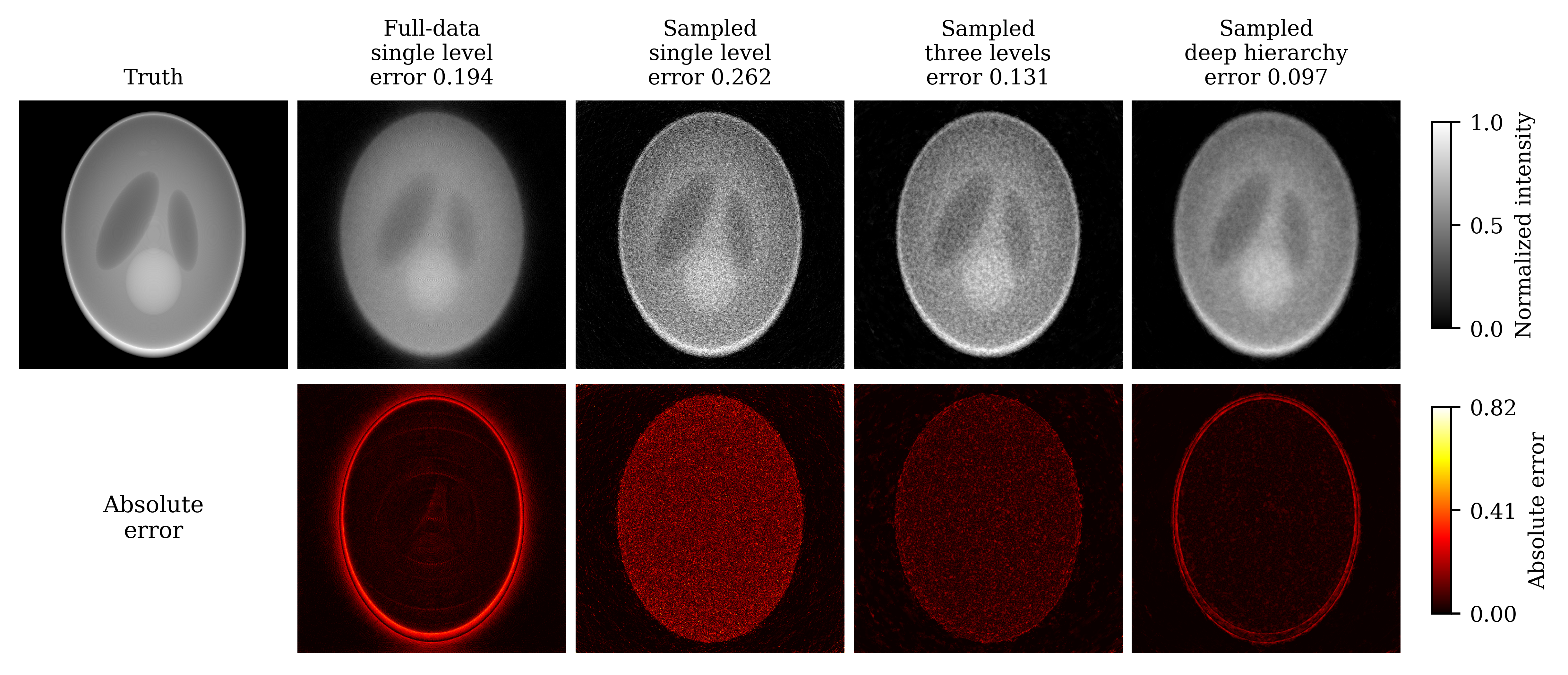}
\caption{Example fixed-angle CT reconstructions at \(N=512\) with \(128\)
projection angles, \(5\%\) relative noise, and LSQR(2) corrections.  The
upper row shows the reconstructions and their relative errors; the lower row
shows the corresponding absolute errors.}
\label{fig:fixed-angle-ct-reconstruction}
\end{figure}

\Cref{fig:fixed-angle-ct-reconstruction} shows example reconstructions for the
\(N=512\) problem with LSQR(2) corrections.  In this example, the full-data
single-level reconstruction is oversmoothed.  Three levels reduce the
artifacts in the sampled single-level reconstruction, and the deep hierarchy
gives the lowest error while preserving the dominant structures.

\paragraph{Matched CT and deblurring comparisons.}
At each of the four batch fractions, both hierarchies improve on the
corresponding single-level kernel in every paired replicate of all \(96\)
comparison groups: two problems, two sizes, two noise levels, three kernels,
two endpoint definitions, and two hierarchy depths.  At the headline
\(0.5\%\) fraction, the deeper hierarchy also improves on three levels in every
pair.  The row and operator endpoints reuse the same trajectories; the six
noise/sampling replicates, not the pooled endpoint count, are the replication
units.

\Cref{tab:large-gradient} displays the gradient control at \(N=511\) and
\(\eta=0.01\).  At the operator-epoch endpoint, the paired single-level minus
six-level error reductions are \(0.07163\,(0.00038)\) for CT and
\(0.03597\,(0.00028)\) for deblurring, with one standard error in parentheses.
Thus the hierarchy improves a direct-gradient baseline even when all inner
operator applications are charged.  Together with the LSQR(2) and LSQR(5)
results, this establishes robustness to the tested local-solve depths.  These
values use the projected Armijo line search in \cref{alg:amls}; in all four cases
the hierarchy advantage remains after the feasible trial point is included in
the line-search test.

\begin{table}[htbp]
{\footnotesize
\caption{Mean relative reconstruction error for gradient corrections at
\(N=511\), \(\eta=0.01\), and \(0.5\%\) iid row batches.  Parentheses give one
standard error over six paired replicates.  CT uses \(511\) angles, and
deblurring uses periodic Gaussian blur with \(\sigma=2\) pixels.  Row and
operator endpoints are the last complete updates within twelve row epochs and
sixty forward/adjoint operator epochs, respectively.}
\label{tab:large-gradient}
\begin{center}
\begin{tabular}{llrrr}
\toprule
Problem & endpoint & single level & three levels & six levels \\
\midrule
CT & row      & 0.1080 (0.0002) & 0.0495 (0.0004) & 0.0357 (0.0004) \\
CT & operator & 0.1074 (0.0002) & 0.0495 (0.0004) & 0.0357 (0.0004) \\
Deblur & row      & 0.0643 (0.0002) & 0.0465 (0.0003) & 0.0241 (0.0001) \\
Deblur & operator & 0.0597 (0.0003) & 0.0439 (0.0003) & 0.0237 (0.0001) \\
\bottomrule
\end{tabular}
\end{center}
}
\end{table}

\begin{figure}[htbp]
\centering
\includegraphics[width=\textwidth]{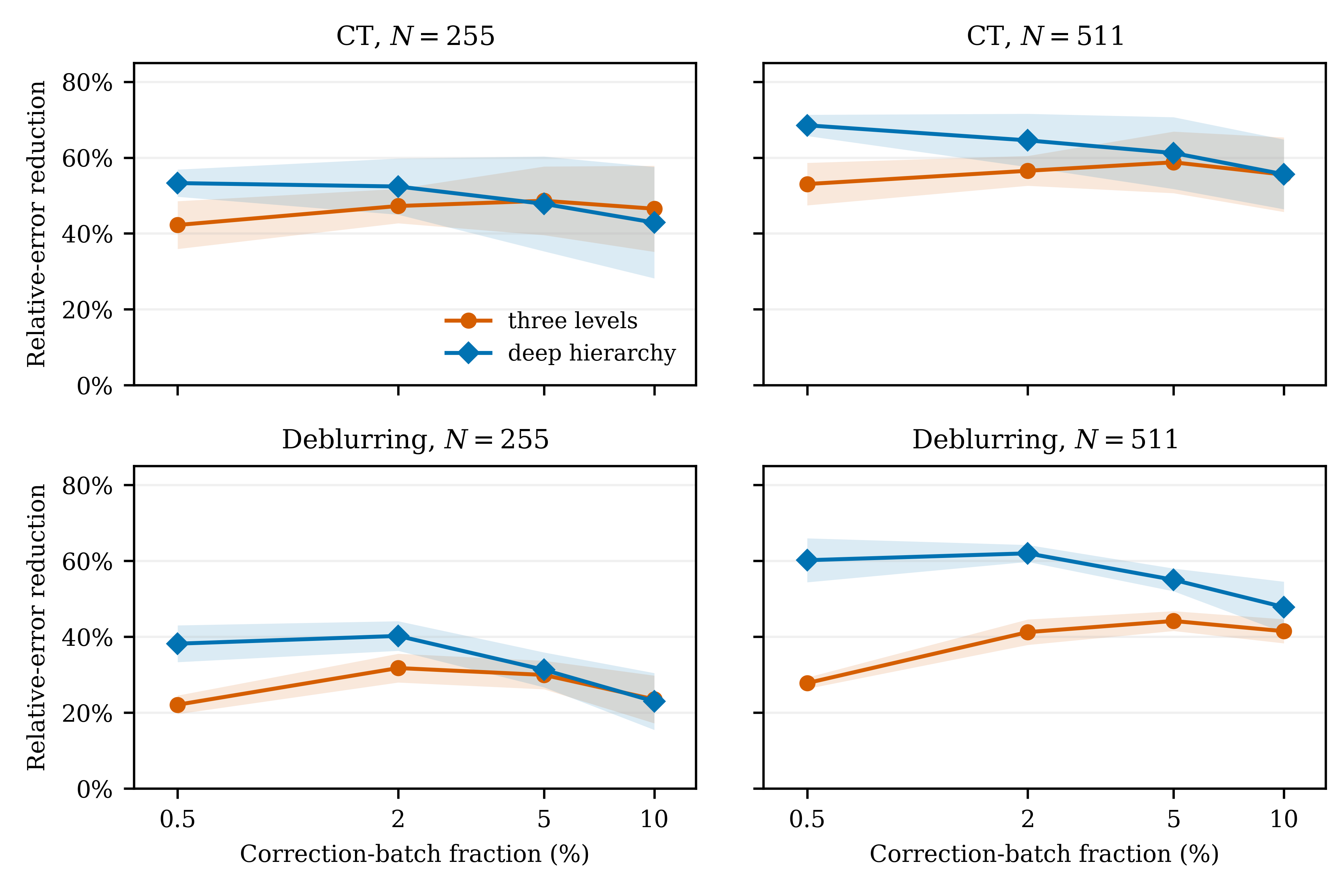}
\caption{Hierarchy gain versus correction-batch fraction at the
sixty-operator-epoch endpoint.  Gain is the percentage reduction in relative
reconstruction error from the matched single-level method.  Markers are means
over two noise levels, three local solvers, and six paired replicates; shaded
bands show one standard deviation across these \(36\) comparisons.}
\label{fig:batch-size-sensitivity}
\end{figure}

The hierarchy gain is not specific to the \(0.5\%\) stress test.
\Cref{fig:batch-size-sensitivity} summarizes the matched operator-epoch
comparisons at all four batch fractions.  The gain remains positive at
\(2\%\), \(5\%\), and \(10\%\) for both problems, grid sizes, hierarchy
depths, noise levels, and local solvers.  The \(2\%\) and \(5\%\) results are
generally comparable to, and sometimes stronger than, those at \(0.5\%\);
at \(10\%\) the gain remains positive but is usually smaller.  The distinction
between the three-level and deeper hierarchies also narrows as the batch grows.
The row-epoch endpoint gives the same qualitative conclusions.

\section{Nonlinear scope study: MAGPIE ptychography}
\label{sec:scope}
\label{subsec:scope-multilevel}

To examine whether hierarchical correction spaces remain useful beyond linear
least-squares problems, we consider known-probe ptychographic phase retrieval.
We compare the MAGPIE surrogate hierarchy~\cite{zhang2026magpie}, motivated by
earlier multigrid optimization for ptychography~\cite{fung2020}, with the same
surrogate update run at a single level.  This is a scope study of the hierarchy,
not a test of the complete stochastic multilevel algorithm: MAGPIE does not use
the sampling strategy or projected Armijo line search developed above.

The tests use an analytic MAGPIE-style surrogate with \(49\) raster scan
positions and \(50\%\) overlap.  For the spectral diagnostic, we use a
\(128\times128\) reference object and a \(32\times32\) probe; the
reconstruction tests use \(256\times256\) and \(512\times512\) objects with
probe width \(N/4\).  These measurements are nonlinear and spatially local:
a scan batch covers part of the object rather than providing exchangeable
global validation data.  The linear theory in \cref{sec:theory} therefore does
not apply directly.

\begin{figure}[htbp]
\centering
\includegraphics[width=\textwidth]{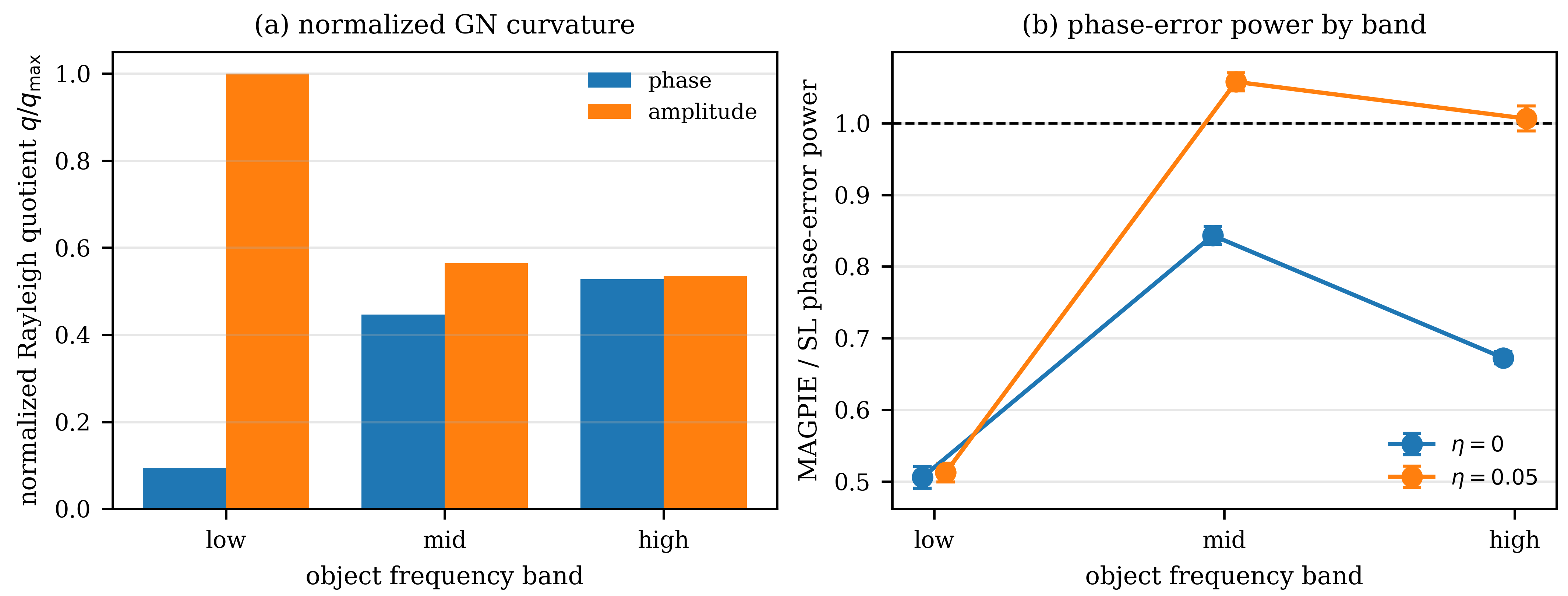}
\caption{Frequency-band diagnostics for the \(128\times128\) MAGPIE reference
problem with a \(32\times32\) probe.  (a) Normalized band-wise curvature of
the object Gauss--Newton Hessian; smaller values indicate more weakly
constrained modes.  (b) Ratio of multilevel to single-level phase-error power;
values below one favor the multilevel surrogate.}
\label{fig:ptycho-mg}
\end{figure}

\Cref{fig:ptycho-mg} identifies the phase modes on which the hierarchy acts.
The band-restricted Rayleigh quotients in panel~(a) show that low-frequency
phase modes are weakly constrained, meaning that perturbations in these modes
produce only small changes in the measured data.  In contrast, low-frequency
magnitude is strongly constrained.
After removing the global phase null mode, the full-batch projected
phase-space condition numbers are \(5.3\), \(65\), and
\(1.46\times10^3\) on \(4\times4\), \(8\times8\), and \(16\times16\)
phase grids, respectively; the coarser phase spaces are therefore better
conditioned.  For an aligned reconstruction \(\widehat x\), we define the
pixelwise phase error as
\(e_\phi=\arg(\widehat x\,\overline{x_\star})\) and its power in a Fourier
frequency band \(\mathcal B\) as
\(P_{\mathcal B}=\sum_{k\in\mathcal B}|\widehat e_\phi(k)|^2\).
Panel~(b) reports the ratio
\(P_{\mathcal B}^{\rm MAGPIE}/P_{\mathcal B}^{\rm single}\), so values below
one indicate less phase error in that band.  The hierarchy gives its largest
reduction in the same weakly constrained low-frequency band.  This agreement supports the
coarse phase-space interpretation, although it does not separate statistical
stabilization from deterministic acceleration or update frequency.

Large-object tests at \(N=256\) and \(512\) provide a separate scope check.
For noiseless intensity \(I\) and \(\eta>0\), the Poisson-noise model
draws \(Z\sim\operatorname{Poisson}(I/\eta)\) and records
\(I_{\eta}=\eta Z\).  Thus the data have mean \(I\) and
variance \(\eta I\), so larger \(\eta\) means stronger noise;
\(\eta=0\) denotes noiseless data.
Across the tested noise levels, scan batches, and both fixed-pass and
fixed-update endpoints, both the five-level hierarchy and the deeper
hierarchy---six levels at \(N=256\) and seven at \(N=512\)---improve on
single level in all three paired seeds.  For \(N=512\),
\(\eta=0.05\), and scan batch \(8\), the paired single-level minus
seven-level error reductions are \(0.0849\) with \(40\) passes and \(0.0724\)
with a fixed update count.  The complete design and representative table are
given in \cref{app:ptycho-transfer}.  These tests support the coarse
phase-space interpretation, but they do not test the sampler or line search
used in the linear experiments.

\section{Discussion}
\label{sec:discussion}

The main finding is that a hierarchy can make sampled reconstruction more
reliable when a local batch contains too little information for a stable
fine-grid update.  This behavior is clearest in the increasingly
ill-conditioned CT refinement: the single-level method degrades as the image
dimension grows, while the multilevel methods maintain or improve their
accuracy.

The benefit is not tied to one local solver or one imaging model.  It persists
for gradient and short Krylov corrections, for CT and periodic deblurring, and
under comparisons based on either sampled rows or forward and adjoint
applications.  Deeper hierarchies generally provide larger gains.  The MAGPIE
study suggests that the same correction-space principle can also be useful in
a nonlinear inverse problem, although it does not test the complete algorithm
developed for the linear setting.

The theory explains this behavior through coarse-space identifiability.  A
useful coarse space represents important components of the reconstruction
error while reducing the number of directions that must be determined from a
batch.  A batch that is insufficient for a stable fine-grid correction may
therefore still determine a reliable coarse correction.  The independent
projected line search then checks the proposed change on separate sampled
data.

These experiments do not establish grid-independent convergence, a wall-clock
speedup, or a universal advantage over full-data optimization.  Instead, they
support a practical principle: use sampled corrections in spaces that are both
expressive enough to capture important errors and small enough to be
identified from the available data.

\appendix
\section{Additional numerical comparisons}
\label{app:additional-results}

\subsection{Correction-space comparison}
\label{app:subspace-comparison}

We compare the geometric space with four spaces of the same dimension:
Haar-random, Fourier/DCT low-pass, leading right singular vectors (RSV) of \(A\),
and \(\mathcal K_{n_\ell}(A^\top A,A^\top b)\).
Let \(Q_V\) be an orthonormal basis for a \(d\)-dimensional test space \(V\),
and set \(B_V^\gamma=Q_V^\top A^\top A Q_V+\gamma I\).  We use the
regularized coherence diagnostic
\begin{equation}
 \mu_V^\gamma
 =\frac{m}{d}\max_{1\leq i\leq m}
 (Q_V^\top a_i)^\top(B_V^\gamma)^{-1}(Q_V^\top a_i).
 \label{eq:regularized-coherence}
\end{equation}
In orthonormal coordinates, we measure approximation defect
\(\|(I-\Pi_V)x^\star\|/\|x^\star\|\), regularized coherence
\(\mu_V^\gamma\), curvature error \(\delta_V^\gamma(S)\), and capture
\(\kappa_V^\gamma(S)\), using \cref{eq:diagnostic-ridge} and identical iid
importance-weighted batches.

\begin{table}[htbp]
{\footnotesize
\caption{Correction-space diagnostics at \(n_\ell=225\) (level~3, \(N=63\)),
averaged over five seeds.  Curvature error uses \(|S|=0.05m\), and ridge-model
capture uses \(|S|=0.02m\); both are evaluated at \(x=0\).  The full-data
Krylov space is a diagnostic reference rather than the sampled Krylov space
used by the reconstruction solver.}
\label{tab:subspace}
\begin{center}
\begin{tabular}{llrrrr}
\toprule
 & Subspace & approx.\ defect & \(\mu_V^\gamma\) & \(\delta_V^\gamma\) & \(\kappa_V^\gamma\) \\
\midrule
CT & geometric & 0.202 & 4.72  & 3.15 & 0.975 \\
   & low-pass (DCT) & 0.206 & 3.63  & 2.99 & 0.976 \\
   & random & 0.970 & 3.69  & 2.70 & 0.646 \\
   & top-RSV & 0.249 & 10.83 & 3.46 & 0.970 \\
   & Krylov & 0.134 & 4.04  & 2.81 & 0.984 \\
\midrule
Deblur & geometric & 0.200 & 2.15 & 6.77 & 0.463 \\
       & low-pass (DCT) & 0.207 & 2.20 & 7.07 & 0.461 \\
       & random & 0.972 & 1.71 & 6.88 & 0.368 \\
       & top-RSV & 0.209 & 1.61 & 6.96 & 0.456 \\
       & Krylov & 0.165 & 1.81 & 6.36 & 0.467 \\
\bottomrule
\end{tabular}
\end{center}
}
\end{table}

The RSV and full-data Krylov spaces are references rather than cost-matched
algorithms or tests of same-batch sampled LSQR
(\cref{rmk:krylov-identifiability}).  The deblurring diagnostic observes a
fixed half of the blurred pixels, unlike the full periodic-blur reconstruction
tests in \cref{sec:kamis}.

\subsection{Large-object MAGPIE reconstruction}
\label{app:ptycho-transfer}

We compare single level, the original five patch levels, and a deeper hierarchy
with six levels at \(N=256\) or seven at \(N=512\).  For a noiseless intensity
\(I\), the Poisson model draws
\(Z\sim\operatorname{Poisson}(I/\eta)\) and records
\(I_{\eta}=\eta Z\), so
\(\mathbb E[I_{\eta}\mid I]=I\) and
\(\operatorname{var}(I_{\eta}\mid I)=\eta I\).
We use \(\eta=0\) for noiseless data.  For each size,
\(\eta\in\{0,0.05\}\), and scan batch \(2\), \(8\), or \(49\), we use
three paired data/sampling seeds.  No independent line search is used.
Relative object error is measured on the central \(N\times N\) region after
complex-scale alignment.

We compare two endpoints.  \(40\) scan passes give \(1000\), \(280\), or \(40\)
updates for batches \(2\), \(8\), or \(49\), respectively.  The fixed-update
control uses \(200\) updates for batches \(2\) and \(49\), and \(203\) for
batch \(8\) because the implementation completes whole passes.  Each
single-level/multilevel comparison uses the same batch and endpoint.

\begin{table}[htbp]
{\footnotesize
\caption{Mean aligned object relative error for known-probe ptychography at
\(N=512\), probe width \(128\), and Poisson scale \(\eta=0.05\).  Means use
three paired data/sampling seeds; all standard errors are below \(10^{-4}\).
The upper block fixes \(40\) scan passes, and the lower block fixes \(200\)
updates (\(203\) at batch \(8\)).}
\label{tab:ptycho-mg-transfer}
\begin{center}
\begin{tabular}{llrrr}
\toprule
Endpoint & batch & single level & five levels & seven levels \\
\midrule
40 passes & 2  & 0.4599 & 0.3908 & 0.3701 \\
          & 8  & 0.4618 & 0.3962 & 0.3769 \\
          & 49 & 0.4716 & 0.4298 & 0.4177 \\
\midrule
200 updates & 2  & 0.4726 & 0.4464 & 0.4368 \\
            & 8  & 0.4652 & 0.4097 & 0.3927 \\
            & 49 & 0.4558 & 0.3723 & 0.3573 \\
\bottomrule
\end{tabular}
\end{center}
}
\end{table}

Batch-size rankings depend on the chosen endpoint; these are matched hierarchy
comparisons, not evidence that smaller scan batches are intrinsically
preferable.

\ifsiamclass
  \bibliographystyle{siamplain}
  \bibliography{references}
\else
  % The desktop editor compiles only this source file, so it cannot read the
  % adjacent BibTeX database.  This generated bibliography mirrors the current
  % references.bib for the isolated preview; the SISC build above remains the
  % canonical BibTeX compilation route.

\fi

\end{document}